\documentstyle[11pt,leqno]{article}
\newtheorem{theorem}{{\sc Theorem}}
\newcommand{\bt}{\begin{theorem}}
\newcommand{\et}{\end{theorem}}
\newcommand{\newsection}[1]{\setcounter{equation}{0} \setcounter{theorem}{0}
\section{#1}}

\newcommand{\bea}{\begin{eqnarray}}
\newcommand{\eea}{\end{eqnarray}}

\def \spec#1 {\mathop{#1}}

\def \b #1 {\bf #1}

\newcommand{\be}{\begin{equation}}
\newcommand{\ee}{\end{equation}}
\newcommand{\ben}{\begin{eqnarray*}}
\newcommand{\een}{\end{eqnarray*}}

\newtheorem{thm}{Theorem}[section]
\newtheorem{cor}[thm]{Corollary}
\newtheorem{lem}[thm]{Lemma}

\newtheorem{rem}[thm]{Remark}

\begin{document}
	\begin{center}\begin {Large}	
A study of the limiting behavior of  delayed random sums under non-identical distributions setup
\end {Large}

M. SREEHARI\\
6-B, Vrundavan Park, New Sama Road, Vadodara, 390024, India
\end{center}
\noindent{\bf Abstract :} We consider delayed sums of the type $S_{n+a_n}-S_n$ where $a_n$ is possibly a positive integer valued random variable satisfying  certain conditions and $S_n$ is the sum of independent random variables $X_n$ with distribution functions $F_n \in \{G_1, G_2\}$. We study the limiting behavior of delayed sums and prove laws of the iterated logarithm of Chover- type. These results extend the results in Vasudeva and Divanji (1992) and Chen (2008).  

\newsection{Introduction and notations} We consider a sequence of independent random variables (rvs) $\{X_n\}$ with corresponding distribution functions $\{F_n\}$ where for each $n$, $F_n \in \{G_1, G_2\}$. We assume that $G_j $ is in the domain of normal attraction of a non-normal stable law with characteristic function $\varphi_j(t)=\exp (-\lambda_j|t|^{\alpha_j})$,  $0 < \alpha_1 \leq \alpha_2 <2$. It is known then that
\begin{equation}\label{E1}
1-G_j(x)=\frac { c_{j1}+\theta_j(x)}{x^{\alpha_j}},\;\;\;G_j(-x)= \frac { c_{j2}+\beta_j(-x)}{x^{\alpha_j}}, \;\;\; x > 0\;
\end{equation}
 where $\theta_j(x),\; \;\beta_j(-x) \rightarrow 0 $ as $x \rightarrow \infty$,  $c_{j1} > 0$,  and $c_{j2}>0$. Set $S_n=\sum_{k=1}^n \;X_k$ and consider the sampling scheme $\{\tau_1(n), \tau_2(n)\}$ where $\tau_j(k)-\tau_j(k-1)=1$ if $F_k=G_j$ and zero otherwise. clearly $\tau_1(n)+\tau_2(n)=n$. Assume that each $\tau_j(n) \rightarrow \infty$. We shall consider the case with  $0 < \alpha_1 < \alpha_2 < 2$ first.  We shall discuss $\alpha_1 = \alpha_2$ case at the end.\\
 
 For later use we introduce the notation $U_{\tau_1(n)}$, the sum of those $X_k$ in $\{X_1, X_2, \cdots , X_n\} $ with distribution function $G_1$ in the domain of normal attraction of the stable $(\alpha_1)$ law and $V_{\tau_2(n)}$, the sum of those $X_k$ in $\{X_1, X_2, \cdots , X_n\} $ with distribution function $G_2$ in the domain of normal attraction of stable  $(\alpha_2)$ law. Then  the limit distribution functions of $ \frac{U_{\tau_1(n)}-d_1(\tau_1(n))}{B_1(\tau_1(n))}$  and $ \frac{V_{\tau_2(n)}-d_2(\tau_2(n))}{B_2(\tau_2(n))}$ are the stable $(\alpha_1)$ and the  stable $(\alpha_2)$ laws respectively for appropriate choices of $d_1(\tau_1(n)) $ and $d_2(\tau_2(n))$. One can choose $d_1(\tau_1(n))=0= d_2(\tau_2(n))$  if $\alpha_1 \neq 1, \alpha_2 \neq 1 $ and  $d_1(n)  \sim  n \;\log n, (d_2(n) \sim  n\;\log n)$  if $\alpha_1 = 1 (\alpha_2 = 1)$ and in the case $\alpha_1=1=\alpha_2 $ we may take $A_n= d_1(\tau_1(n))+\;d_2(\tau_2(n))$.  Here we follow the notation $ f_n \;\sim \; g_n$  if $f_n/g_n \rightarrow C$, as $ n \rightarrow \infty$, where $ 0\;<\;C\;<\infty$.  Henceforth we assume that the limit distribution of $\frac{(S_n-A_n)}{B_n}$ exists. Thus if the limit distribution is a convolution of the two stable laws  then $0\;<\;\alpha_1\;<\;\alpha_2\;<\;2$ and $\tau_1(n) \sim n^{\alpha_2/\alpha_1}$ and $ \tau_2(n)\sim n$. If the limit distribution is stable ($\alpha_2$), $\tau_2(n)\sim n$ because $\frac{(\tau_1(n))^{\alpha_2}}{(\tau_2(n))^{\alpha_1}}=\left( \frac{\tau_1(n)} {\tau_2(n)}\right)^{\alpha_1}[\tau_1(n)]^{\alpha_2-\alpha_1}\rightarrow 0 $. If the limit distribution is stable ($\alpha_1$) then $n\;(\tau_1(n))^{-\;\alpha_2/\alpha_1} \rightarrow 0$.  Unfortunately  no more specific behavior can be made out about $\tau_j$s when the limit distribution is stable ($\alpha_1$) or stable ($\alpha_2$)  as in the case where the limit distribution is a composition.  We therefore need suitable condition in the case the limit distribution is stable ($\alpha_1$).\\\\
 Chover (1966) was the first to prove a law of the iterated logarithm (LIL) for the symmetric stable laws with exponent $\alpha < 2$ where he considered the limiting behavior of $ \left| \frac{S_n}{n^{1/\alpha}}\right| ^{1/\log \log n}$ and Heyde (1969) extended Chover's result to certain rvs with common distribution in the domain of normal attraction of the symmetric stable law with  exponent $\alpha \neq 1, 2.$  Zinchenko (1994)  extended Chover's LIL  for independent identically distributed (iid) symmetric stable ($0\;<\;\alpha\;<\;2$) rvs. \\\\
 Consider the delayed sums $T_n=S_{n+a_n} \;- S_n$ where $\{ a_n\rightarrow \infty\} $ is a sequence of positive integers. Lai (1974) proved the LIL for delayed sums.  Vasudeva and Divanji (1993) extended the result of Chover to the non-identical distribution setup assuming $G_j, j=1, 2$ to be positive stable laws with exponents $ 0 < \alpha_1 \leq \alpha_2 < 1$. They assumed that the limit distribution of $S_n$, properly normed, exists and is a composition of the two stable laws. Chen (2008)   proved some  general  results on the limiting behavior of $S_n$ and derived extension of the result of   Vasudeva and Divanji (1993) to the case  of symmetric stable laws $G_j$ with exponents $0 < \alpha_1 \leq \alpha_2 < 2$ thereby relaxing the assumption of positive stable laws   $G_j, j=1, 2$. Henceforth we drop the term symmetric and just refer to the limit distributions as stable laws. \\\\
 
 The main aims of this paper are:\\
 (i) to extend the results to the case where each $F_n$ is in the domain of normal attraction of the  stable law $G_1$ or $G_2$ according to the sampling scheme described above and satisfying certain conditions. We shall not restrict to the case of the limit distribution of $S_n$, properly normed, is a composition of the two  stable laws; that is, the limit distribution may be  stable ($\alpha_1$) or  stable ($\alpha_2$),\\and \\
 (ii) to extend the results of Chen (2008) to the case where the lags $a_n$ are positive rvs independent of the summands $X_k$ in the context described in (i).  \\

  In Section 2 we  state the results of Chen (2008).  Further we prove an extension of Lemma 2.1 in Chen (2002). In Section 3 we  discuss the  delayed sum problem  when  $F_n$ are in the domains of attraction of stable laws and in Section 4 we consider  similar problems with random $a_n$.
 
\newsection{Statements of Chen's results} Chen (2008) investigated the almost sure limiting behavior of partial sums $S_n $ and proved Chover's LIL type results for  the delayed sums $T_n$ under the assumption that $G_j$ are non-normal stable. For the sampling scheme $\left\lbrace \tau_1(n), \tau_2(n) \right\rbrace $ a necessary and sufficient condition for $\frac{(S_n-A_n)}{B_n}$, with $A_n \in R$ and $B_n > 0$,  to converge in distribution to a proper rv is that the ratio  $\frac{(\tau_1(n))^{\alpha_2}}{(\tau_2(n))^{\alpha_1}}\rightarrow \lambda$, where $\lambda \geq 0$.  If $\lambda= 0$  the limit distribution is  the stable $(\alpha_2)$,  if $0 < \lambda < \infty$ the limit distribution is a composition of the stable laws with exponents $\alpha_1$ and $\alpha_2$ and if $\lambda =\infty$ the limit distribution is the stable $(\alpha_1)$. In the   case  of $\lambda=\;\infty$ we may take
$B_n \sim B_1(\tau_1(n)) \sim (\tau_1(n))^{1/\alpha_1}. $   In the case the limit distribution is the  stable  $(\alpha_2) $    we may take   $B_n \sim B_2(\tau_2(n))=(\tau_2(n))^{1/\alpha_2} \sim n^{1/\alpha_2}$. Further when the limit distribution is a composition of the two stable laws $\tau_1(n)\sim [n^{\alpha_1/\alpha_2}],\; \tau_2(n)\sim n $ and we may take $B_n\sim  n^{1/\alpha_2}.$  For details we refer to Sreehari (1970).\\\\ 
We now introduce some assumptions  which are assumed in different situations: \\
Assumption ($C_1$): $limsup_{n\rightarrow \infty}\; a_n/\tau_1(n)\;<\,\infty.$\\
Assumption ($C_2$) $limsup_{n\rightarrow \infty}\; a_n/n\;<\,\infty$.\\
Assumption ($C_3$) For some $\mu > \frac{\alpha_2-\alpha_1}{\alpha_2},\;\;\;\tau_1(n) < n^{\alpha_1/\alpha_2}\;(\log n)^{-\mu}.$\\
  Note that  the assumption that $limsup_{n\rightarrow \infty} \;a_n/\tau_1(n)\;<\,\infty$ is slightly stronger than the assumption  $limsup_{n\rightarrow \infty} \;a_n/n\;<\,\infty$ which was assumed by Chen (2008) in the case $0\;<\;\lambda\;<\;\infty.$ We shall   assume ($C_2$) while dealing with the case $0 \leq \lambda < \infty $ and ($C_1$) while dealing with the case $\lambda = \infty$.\\\
  We shall now recall Chen's results who assumed, like Vasudeva and Divanji,  that the above limit distribution is a composition of the two stable  laws with exponents $\alpha_1$ and $\alpha_2$ and proved the following.
\begin{thm}(Chen, 2008) \label{T1} Let $f > 0$ be a  nondecreasing function. Then with probability one
\[\limsup_{n\rightarrow \infty}\frac{|S_n|}{B_n\;(f(n))^{1/\alpha_1}}=\left\{\begin{array}{ll}
0 & \mbox{}\\
\infty & \mbox{}
\end{array}\right.
\;\;\; \Longleftrightarrow \;\;\; \int_1^\infty \frac{1}{xf(x)}dx \;\;\;\left\{\begin{array}{ll}
<\infty & \mbox{}\\
=\infty. & \mbox {}
\end{array}\right.\]

\end{thm}

\begin{cor} \label{C1} For  every $\delta > 0$ 
	$$\limsup_{n\rightarrow \infty}\frac{|S_n|}{B_n(\log n)^{(1+\delta)/\alpha_1}}= 0 \;\;\;a.s. $$
	and $$\limsup_{n\rightarrow \infty}\frac{|S_n|}{B_n(\log n)^{1/\alpha_1}}=\infty \;\;\;a.s. $$
	In particular $$\lim \sup_{n\rightarrow \infty}\left| \frac{S_n}{B_n}\right| ^{1/\log \log n}= e^{1/\alpha_1} \;\;\;a.s. $$
	\end{cor}
\begin{rem}
(1)	When the limit distribution of $S_n$ is stable ($\alpha_1$)  also the same proof of Theorem \ref{T1} and Corollary \ref{C1} will go through with minor modifications.\\
(2) When the limit distribution is stable ($\alpha_2$), under  the Assumptions ($C_2$) and ($C_3$) the following result holds. 
For  every $\delta > 0$ 
$$\limsup_{n\rightarrow \infty}\frac{|S_n|}{B_n(\log n)^{(1+\delta)/\alpha_2}}= 0 \;\;\;a.s. $$
and $$\limsup_{n\rightarrow \infty}\frac{|S_n|}{B_n(\log n)^{1/\alpha_2}}=\infty \;\;\;a.s. $$
In particular $$\lim \sup_{n\rightarrow \infty}\left| \frac{S_n}{B_n}\right| ^{1/\log \log n}= e^{1/\alpha_2} \;\;\;a.s. $$
 
\end{rem}
\begin{thm}(Chen, 2008) \label{T2} Let $f > 0$ be a nondecreasing function and let $\{a_n\} $ satisfy  the Assumption ($C_2$). Then with probability one
	\[\limsup_{n\rightarrow \infty}\frac{|T_n|}{B_n\;(f(n))^{1/\alpha_1}}=\left\{\begin{array}{ll}
	0 & \mbox{}\\
	\infty & \mbox{}
	\end{array}\right.
	\;\;\; \Longleftrightarrow \;\;\; \int_1^\infty \frac{1}{xf(x)}dx \;\;\;\left\{\begin{array}{ll}
	<\infty & \mbox{}\\
	=\infty & \mbox {}
	\end{array}\right.\]
	
\end{thm}

\begin{cor} \label{C2} Let $\{a_n\}$ satisfy the Assumption ($C_2$). Then for  every $\delta > 0$ 
	$$\limsup_{n\rightarrow \infty}\frac{|T_n|}{B_n(\log n)^{(1+\delta)/\alpha_1}}= 0 \;\;\;a.s. $$
	and $$\limsup_{n\rightarrow \infty}\frac{|T_n|}{B_n(\log n)^{1/\alpha_1}}=\infty \;\;\;a.s. $$
	In particular $$\lim \sup_{n\rightarrow \infty}\left| \frac{T_n}{B_n}\right| ^{1/\log \log n}= e^{1/\alpha_1} \;\;\;a.s. $$
\end{cor}
\begin{rem} 
(1)	In the case that the limit distribution of $S_n$, properly normed, is  the stable $(\alpha_1)$ law   we may take $B_n=(\tau_1(n))^{1/\alpha_1}.$ The same results hold if the  Assumption ($C_1$) holds.\\
(2) In the case that the limit distribution of $S_n$, properly normed, is  the stable $(\alpha_2)$ law  we may take $B_n=(\tau_2(n))^{1/\alpha_2}.$ Then the following result holds:\\
Let $\{a_n\}$ satisfy the Assumptions ($C_2$) and ($C_3$). Then for  every $\delta > 0$ 
$$\limsup_{n\rightarrow \infty}\frac{|T_n|}{B_n(\log n)^{(1+\delta)/\alpha_2}}= 0 \;\;\;a.s. $$
and $$\limsup_{n\rightarrow \infty}\frac{|T_n|}{B_n(\log n)^{1/\alpha_2}}=\infty \;\;\;a.s. $$
In particular $$\lim \sup_{n\rightarrow \infty}\left| \frac{T_n}{B_n}\right| ^{1/\log \log n}= e^{1/\alpha_2} \;\;\;a.s. $$

\end{rem}
\begin{cor} \label{C3} Let $\{a_n\}$ be a subsequence of positive integers with $\limsup_{n\rightarrow \infty}a_n/n < \infty$ and let $\gamma_n=\log(n/a_n)+ \log \log n$. \\
	(i) If $\lim_{n \rightarrow \infty}\frac{\log (n/a_n)}{\log \log n}=\infty$, then
	$$\limsup_{n\rightarrow \infty}\left| \frac{T_n}{B_{a_n}}\right| ^{1/\gamma_n}=e^{1/\alpha_2} \;\;\; a.s.$$
	(ii) If $\lim_{n \rightarrow \infty}\frac{\log (n/a_n)}{\log \log n}=0$, then
	$$\limsup_{n\rightarrow \infty}\left| \frac{T_n}{B_{a_n}}\right| ^{1/\gamma_n}=e^{1/\alpha_1} \;\;\; a.s.$$
	(iii) If $\lim_{n \rightarrow \infty}\frac{\log (n/a_n)}{\log \log n}=s \in (0, \infty)$, then
	$$\limsup_{n\rightarrow \infty}\left| \frac{T_n}{B_{a_n}}\right| ^{1/\gamma_n}=e^{\frac{\alpha_1 s+\alpha_2}{(s+1)\alpha_1\alpha_2}} \;\;\; a.s.$$
\end{cor}

The proofs of these results heavily depend on the fact that  $ \frac{U_{\tau_1(n)}-b_{\tau_1(n)}}{(\tau_1(n))^{1/\alpha_1}}$  and $ \frac{V_{\tau_2(n)}-d_{\tau_2(n)}}{(\tau_2(n))^{1/\alpha_2}}$ are distribute as stable ($\alpha_1$) and stable ($\alpha_2$) respectively. This does not hold in the case $G_j$ is not stable as in our case. To circumvent this difficulty we use the  lemma \ref{L1} below. In the rest of the paper we denote $C$ as a generic positive number which may be different at different places. Before we close this Section we shall prove an extension of the result in Lemma 2.1 in Chen (2002) for a sequence of independent rvs $ \{Z_k\}$ with the common distribution function H in the domain of normal attraction of the stable law with  characteristic function $\varphi(t)=\exp (-\lambda|t|^{\alpha})$. We denote $W_n=Z_1+Z_2+\cdots+Z_n$. Then we have the following.
\begin{lem} \label{L1}
	Let $f > 0$ be a nondecreasing function satisfying $\int_1^\infty \frac{1}{x\;f(x)} dx < \infty .$  Then $$\lim_{n\rightarrow \infty} \frac{max_{1\leq k \leq n}|W_k|}{(n\;f(n))^{1/\alpha}}=\;\;0\;\;\; a.s. $$
	\end{lem} 
\bf{Proof.} For any $\epsilon > 0$, let\\ $E_n=\{max_{1 \leq k \leq n}|W_k |> \epsilon\;(nf(n))^{1/\alpha}\}$  and $E_n^*=\{max_{2^n \leq k < 2^{n+1}}\;|W_k |> \epsilon\;(2^nf(2^n))^{1/\alpha}\}$.\\ Then $\limsup_{n \rightarrow \infty} E_n\subset \limsup_{n \rightarrow \infty} E_n^*. $
By the L\'{e}vy inequality, we have for all $n \geq 1$, $P(E_n^*) \leq 2 P(D_n)$ where $D_n= \{|W_{2^{n+1}-1}| \;> \epsilon\;(2^n\;f(2^n))^{1/\alpha} \}$. Since $Z_k$ follows $H$ and $H$ is of the same type as $G$ in (\ref{E1}) we have $$P(D_n) = (2^{n+1}-1) \frac {C +\theta(\epsilon\;(2^n\;f(2^n))^{1/\alpha})+\beta(-\epsilon\;(2^n\;f(2^n))^{1/\alpha})}{\epsilon^\alpha\;2^n\;f(2^n) }.$$ Hence for $N$ sufficiently large 
		$$\sum_{n=N}^\infty P(D_n) < \sum_{n=N}^\infty \frac{C}{f(2^n) } < \int_1^\infty \frac{1}{x\;f(x)} dx <  \infty. $$
\newsection{New results for delayed sums}  We assume that the independent  rvs  $\{X_n\}$ have corresponding distribution functions $\{F_n\}$  where for each $n$, $F_n \in \{G_1, G_2\}$. In the following Lemma we assume that $G_j $ is in the domain of  attraction of a  stable law with characteristic function $\varphi_j(t)=\exp (-\lambda_j|t|^{\alpha_j})$,  $0 < \alpha_1 < \alpha_2 <2$.  Then we have the following  
\begin{lem}\label{L2}
	For any positive constant $ M$  and non-decreasing function $f > 0$ if  
	 $$\int_1^\infty \frac{1}{x\;f(x)} dx = \infty $$
then $$\sum_{n=1}^\infty P(|X_n| \geq M\;B_n(f(n))^{1/\alpha_1}) = \infty. $$
\end{lem} 
\bf{Proof.} Recall that $(\tau_1(n))^{\alpha_2} / (\tau_2(n))^{\alpha_1}\rightarrow \lambda.$  We consider the case with the  $0 \leq \lambda < \infty$ first  and  in this case we may take $B_n= B_2(\tau_2(n))=\;(\tau_2(n))^{1/\alpha_2}$.  We recall that $X_n$ follows $ G_2$ if $\tau_2(n) - \tau_2(n-1)= 1$. Then we have
\begin{eqnarray*}
	\sum_{n=1}^\infty P(|X_n| \geq M\;B_n(f(n))^{1/\alpha_1})& = & \sum_{k=0}^\infty\; \sum_{n=2^k}^{2^{k+1}-1}P(|X_n| \geq M\;B_n(f(n))^{1/\alpha_1})\\
	& \geq & \sum_{k=K_0}^\infty\; \sum_{n=2^k}^{2^{k+1}-1}\frac{C\;L_2(M\;B_n(f(n))^{1/\alpha_2})}{B_n^{\alpha_2}(f(n))^{\alpha_2/\alpha_1}}\\
		& \geq & C\; \sum_{k=K_0}^\infty \sum_{n=2^k}^{2^{k+1}-1} (f(n))^{(\alpha_2-\theta)/\alpha_1}\frac{L_2(B_n)}{B_n^{\alpha_2} \;(f(n))^{\alpha_2/\alpha_1}}
		\end{eqnarray*}\\
where $L_2$ is a slowly varying function,   $K_0$ large and   $\theta > 0$ small by Potter's inequality for   regularly varying functions. ( See  Proposition B.1.9(5), p. 367, De Haan and Ferreira, 2006). The penultimate inequality is obtained by omitting the terms that involve the rvs $X_n$ that follow $G_1$ and then using the well-known relation (8.6) on page 313 in Feller (1970) for the tail probability of the distributions attracted to the stable ($\alpha_2$) law. Then using the fact that as $n\rightarrow \infty$ $\frac{n\;L_2(B_2(n))}{(B_2(n))^{\alpha_2}}\rightarrow C\;>\;0$ and recalling that $B_n=B_2(\tau_2(n))$ we have 
\begin{eqnarray*} 
\sum_{n=1}^\infty P(|X_n| \geq M\;B_n(f(n))^{1/\alpha_1})	&\geq & C\; \sum_{k=K_0}^\infty \;\sum_{n=2^k}^{2^{k+1}-1} \frac{\tau_2(n)\;L_2(B_2(\tau_2(n)))}{(B_2(\tau_2(n)))^{\alpha_2}}\frac{1}{\tau_2(n)(f(n))^{\theta/\alpha_1}}\\
		& \geq & C\; \sum_{k=K_1}^\infty\; \sum_{n=2^k}^{2^{k+1}-1} \frac{1}{\tau_2(n)(f(n))^{\theta/\alpha_1}}\\
		\end{eqnarray*}
	\begin{eqnarray*}	
		& \geq & C \; \sum_{k=K_1}^\infty\;\left[ \tau_2(2^{k+1}-1)-\tau_2(2^k-1)\right] \frac{1}{\tau_2(2^{k+1})\;f(2^{k+1})}\\
\end{eqnarray*}
	\begin{eqnarray}\label{E2}
	& \geq & C \;\sum_{k=K_1+1}^\infty\;\frac{1}{f(2^{k})}
			\end{eqnarray}
for $K_1 > K_0$ since $\theta $ can be chosen to be $< \alpha_2$.  \\
Next we note that
\begin{eqnarray*} 
\int_1^\infty \frac{1}{x\;f(x)} dx & = & \sum_{k=0}^\infty\int_{x=2^k}^{2^{k+1}-1}  \frac{1}{x\;f(x)} dx\\
& \leq & C\;\sum_{k=0}^\infty \frac{1}{f(2^k)} \leq C\;\sum_{k=K_1+1}^\infty \frac{1}{f(2^k)}.\\
\end{eqnarray*}
This together with (\ref{E2}) completes the proof of the Lemma in the case $0 \leq \lambda < \infty$. Steps in the  case of $\lambda = \infty$  can be written on the same lines by recalling that $B_n =  B_1(\tau_1(n))$ and considering the terms for which $X_n$ follows $G_2$ in the summation in stead of those for which $X_n$ follows $G_1$ while deriving the inequality (\ref{E2}). \\
\begin{rem} We recall that this Lemma is proved under the assumption that  $G_j$ is in the domain of attraction of the stable law ($\alpha_j$) and hence in a more general set up than  for the other results.
	\end{rem}
Our next result shows that Theorem $\ref{T1}$ holds when $G_j$ is in the domain of normal attraction of the stable ($\alpha_j$) law for $j=1, 2$ with $0 < \alpha_1 < \alpha_2 < 2.$ 
\begin{thm} \label{T4} Let $f > 0$ be a  nondecreasing function and let $0 < \lambda \leq \infty$. Then with probability one
	\[\limsup_{n\rightarrow \infty}\frac{|S_n|}{B_n\;(f(n))^{1/\alpha_1}}=\left\{\begin{array}{ll}
	0 & \mbox{}\\
	\infty & \mbox{}
	\end{array}\right.
	\;\;\; \Longleftrightarrow \;\;\; \int_1^\infty \frac{1}{xf(x)}dx \;\;\;\left\{\begin{array}{ll}
	<\infty & \mbox{}\\
	=\infty. & \mbox {}
	\end{array}\right.\]
In the case $\lambda = 0$, i.e., when  the limit distribution is stable ($\alpha_2$), if the Assumption ($C_3$) holds we have with probability one 
\[\limsup_{n\rightarrow \infty}\frac{|S_n|}{B_n\;(f(n))^{1/\alpha_2}}=\left\{\begin{array}{ll}
0 & \mbox{}\\
\infty & \mbox{}
\end{array}\right.
\;\;\; \Longleftrightarrow \;\;\; \int_1^\infty \frac{1}{xf(x)}dx \;\;\;\left\{\begin{array}{ll}
<\infty & \mbox{}\\
=\infty. & \mbox {}
\end{array}\right.\]	
\end{thm}
\bf{Proof.} In the following steps $B_n=B_1(\tau_1(n))= (\tau_1(n))^{1/\alpha_1}. $   Assume that $\int_1^\infty \frac{1}{x\;f(x)}dx < \infty$. Clearly $\frac{\log n}{f(n)} \rightarrow 0$ and hence $f(n) \rightarrow \infty$ as $n \rightarrow \infty$. By symmetrization argument (see Lemma 3.2.1 in Stout, 1974) we can prove the result assuming $X_n$s to be symmetric.
 Now by Lemma \ref{L1}
\begin{equation} \label{E3}
 \limsup_{n \rightarrow \infty}\frac{|U_{\tau_1(n)}|}{(\tau_1(n)\;f(\tau_1(n)))^{1/\alpha_1}}  \leq  \limsup_{n \rightarrow \infty}\frac{\max_{1 \leq k \leq \tau_1(n)} |U_k|}{(\tau_1(n)\;f(\tau_1(n)))^{1/\alpha_1}} = 0\;\;\;a.s.\\
  \end{equation}
 Similarly, 
\begin{equation} \label{E4}
\limsup_{n \rightarrow \infty}\frac{|V_{\tau_2(n)}|}{(\tau_2(n)\;f(\tau_2(n)))^{1/\alpha_2}}= 0 \;\;\; a.s.
\end{equation}
 Hence proceeding as in Chen (2008) we have from $ (\ref{E3})$ and $(\ref{E4})$ and the facts $f(\tau_j(n))\leq f(n)$ and $0\;<\alpha_1\;<\;\alpha_2\;<\;2$
\begin{eqnarray*} 
\limsup_{n \rightarrow \infty}\frac{|S_n|}{B_n(f(n))^{1/\alpha_1}} & \leq & \limsup_{n \rightarrow \infty}\frac{|U_{\tau_1(n)}|}{B_n(f(n))^{1/\alpha_1}}\;+\;\limsup_{n \rightarrow \infty}\frac{|V_{\tau_1(n)}|}{B_n(f(n))^{1/\alpha_1}} \\
& \leq & \limsup_{n\rightarrow \infty}\frac{(\tau_1(n)f(\tau_1(n)))^{1/\alpha_1}}{B_n(f(n))^{1/\alpha_1}}\frac{|U_{\tau_1(n)}|}{(\tau_1(n)f(\tau_1(n)))^{1/\alpha_1}} \\
& + & \limsup_{n\rightarrow \infty}\frac{(\tau_2(n)f(\tau_2(n)))^{1/\alpha_2}}{B_n(f(n))^{1/\alpha_1}}\frac{|V_{\tau_2(n)}|}{(\tau_2(n)f(\tau_2(n)))^{1/\alpha_2}}\\
& = & 0\;\;\;a.s.\\
\end{eqnarray*}
Here we use the facts that if $\lambda= \infty $, 
$\frac{(\tau_1(n))^{1/\alpha_1}}{B_n}= 1 $ and $\frac{(\tau_2(n))^{1/\alpha_2}}{B_n}\rightarrow 0$ and if $0 < \lambda < \infty$, $\frac{(\tau_1(n))^{1/\alpha_1}}{B_n}\rightarrow\lambda^{1/\alpha_1\alpha_2} $ while $\frac{(\tau_2(n))^{1/\alpha_2}}{B_n}= 1 $.  \\
We now turn to the divergence part. Assume that $\int_1^\infty \frac{1}{x\;f(x)}dx = \infty$. By lemma \ref{L2} we then have for any $ M > 0$
 \begin{equation}\label{E5}
\sum_{n=1}^\infty P(|X_n| \geq M\;B_n(f(n))^{1/\alpha_1}) = \infty
\end{equation}
 which by Borel -Cantelli lemma implies
 \begin{equation}
 \limsup_{n\rightarrow \infty}\frac{|X_n|}{B_n(f(n))^{1/\alpha_1}} = \;\infty \;\;\;a.s.\
 \end{equation}
Note that  
\begin{eqnarray*}
\limsup_{n\rightarrow \infty}\frac{|X_n|}{B_n(f(n))^{1/\alpha_1}} & \leq  & \limsup_{n\rightarrow \infty}\frac{|S_n|}{B_n(f(n))^{1/\alpha_1}}\\
& +& \limsup_{n\rightarrow \infty}\frac{B_{n-1}\;(f(n-1))^{1/\alpha_1}}{B_n\;(f(n))^{1/\alpha_1}} \frac{|S_{n-1}|}{B_{n-1}\;(f(n-1))^{1/\alpha_1}} \\
& \leq & 2 \limsup_{n\rightarrow \infty}\frac{|S_n|}{B_n(f(n))^{1/\alpha_1}}\\
\end{eqnarray*} and hence from (3.4) we have $$\limsup_{n\rightarrow \infty}\frac{|S_n|}{B_n(f(n))^{1/\alpha_1}}\;= \infty \;\;\; a.s.$$
In the case $\lambda = 0$ similar steps give the result with $\alpha_1$ replaced by 
$\alpha_2$ at appropriate places and using the Assumption ($C_3$).
\begin{rem}  For $f(x)=\log x$, $$\int_1^\infty \frac{1}{x(f(x))^\eta}dx $$  is finite or infinite according as $\eta > 1$ or $\leq 1$. Hence by Lemma 3.1 in Li and Chen (2014) we note that Corollary \ref{C1} will hold in the case $G_j, j=1, 2$ is in the domains of normal attraction of the corresponding stable laws. Thus Theorem \ref{T1} and  corollary \ref{C1} follow from  the above Theorem. Furthermore, Corollary \ref{C1} holds with $\alpha_1$ replaced by $\alpha_2$ in the case $\lambda=0$ when $G_j,\;j=1, 2$, are in the domains of attraction of the corresponding stable laws.\\
\end{rem}
	We now give an extension of Theorem \ref{T2} to the situation where $G_j, j=1, 2$ are in the domains of normal attraction of the stable laws with characteristic functions $\varphi_j(t)=\exp (-\lambda_j|t|^{\alpha_j})$,  $0 < \alpha_1 < \alpha_2 <2, j=1, 2$ . \\
\begin{thm} \label{T5}
Let $f > 0$ be a nondecreasing function such that $limsup_{n\rightarrow \infty}f(b_n)/f(n) < \infty$ if $limsup_{n\rightarrow \infty} b_n/n < \infty$  and let $\{a_n\} $ be a subsequence of positive integers satisfying the Assumption ($C_1$) in the case $\lambda = \infty$, the Assumption ($C_2$) in the case $ 0 < \lambda < \infty$. Then with probability one
\[\limsup_{n\rightarrow \infty}\frac{|T_n|}{B_n\;(f(n))^{1/\alpha_1}}=\left\{\begin{array}{ll}
0 & \mbox{}\\
\infty & \mbox{}
\end{array}\right.
\;\;\; \Longleftrightarrow \;\;\; \int_1^\infty \frac{1}{xf(x)}dx \;\;\;\left\{\begin{array}{ll}
<\infty & \mbox{}\\
=\infty. & \mbox {}
\end{array}\right.\]
 In the case $\lambda = 0$ let $\{a_n\} $ be a subsequence of positive integers satisfying the Assumptions ($C_2$) and  ($C_3$). Then with probability one 
 \[\limsup_{n\rightarrow \infty}\frac{|T_n|}{B_n\;(f(n))^{1/\alpha_2}}=\left\{\begin{array}{ll}
 0 & \mbox{}\\
 \infty & \mbox{}
 \end{array}\right.
 \;\;\; \Longleftrightarrow \;\;\; \int_1^\infty \frac{1}{xf(x)}dx \;\;\;\left\{\begin{array}{ll}
 <\infty & \mbox{}\\
 =\infty. & \mbox {}
 \end{array}\right.\]
	\end{thm}
We omit the proof as it is exactly on the same lines as in Chen (2008) and by using Lemma \ref{L2} in the divergence part. Further, it is a particular case of Theorem \ref{T7} proved in the next Section.
\begin{rem} Corollaries \ref{C2} and \ref{C3} can be easily deduced under the weaker assumption that $G_j$ is in the domain of normal attraction of the stable ($\alpha_j$) law when the limit distribution of $S_n$, properly normed, is a composition of the two stable laws.
	\end{rem}
\newsection{Delayed random sums}
We shall now  consider the situation where each $a_n $ may be a positive integer valued rv. Very  little work  is done in this set up. There is however a large body of work related to  random sums and  random indexed statistics. The importance of this area of research is seen in reliability, insurance, financial mathematics and statistical quality control. We envisage that the delayed random sum theory will have applications in studies concerning control charts with censored samples where the sample size on each occasion will be a random number. To the best of our knowledge there are only two papers dealing with this kind of problem, viz., Divanji and Raviprkash (2016) and Divanji (2017). Both the papers  deal with positive valued rvs which are  identically distributed under rather strange assumptions/conditions. The usual method of investigation in limit theorems with random index is to convert them to limit theorems for non-random index and apply existing results. This is usually  done via what is known as Anscombe's condition or Gnedenko's Transfer theorem. But these techniques of conversion from random index to non-random index do not seem to work  in almost sure limit theory except when the original random variables are positive valued. However, the method proposed by Chen (2008) helps us dealing with random index in LIL discussed in this Section. We impose slightly stronger conditions on the rvs $a_n$ than those in  Theorem \ref{T5}. Our first result below is a direct application of Theorem 2.2 in Gut (2009).\\ Let us introduce the following assumptions:\\
Assumption ($C_1^*$): $limsup_{n\rightarrow \infty}\; a_n/\tau_1(n)\;<\,\infty \;\;\; a.s.$\\
Assumption ($C_2^*$) $limsup_{n\rightarrow \infty}\; a_n/n\;<\,\infty \;\;\; a.s.$\\

\begin{thm} \label{T6}
Let $f > 0$ be a  nondecreasing  function such that $limsup_{n\rightarrow \infty}f(b_n)/f(n) < \infty$ if $limsup_{n\rightarrow \infty} b_n/n < \infty$  and let $ \{a_n\}$ be a sequence of positive integer valued rvs, independent of the rvs $X_k$. Then under the Assumption ($C_2^*$) if $0 < \lambda < \infty$ and the Assumptions ($C_1^*$)  if $ \lambda = \infty$ we have with probability one
$$\limsup_{n\rightarrow \infty}\frac{|S_{n+a_n}|}{B_n\;(f(n))^{1/\alpha_1}}=0 \;\;\; or \; 1$$
according as   $$\int_1^\infty \frac{1}{x f(x)}dx < \infty \;\;\; or \;= \infty. $$
Further if $ \lambda = 0$, under assumptions ($C_2^*$) and  ($C_3$),  we have with probability one
$$\limsup_{n\rightarrow \infty}\frac{|S_{n+a_n}|}{B_n\;(f(n))^{1/\alpha_2}}=0 \;\;\; or \; 1$$
according as  $$\int_1^\infty \frac{1}{x f(x)}dx < \infty \;\;\;or \;= \infty. $$
 
	\end{thm}
\bf{Proof.} Recall from Theorem \ref{T1} that $$\limsup_{n\rightarrow \infty}\frac{|S_{n}|}{B_n\;(f(n))^{1/\alpha_1}}=0 $$ if $$\int_1^\infty \frac{1}{x f(x)}dx < \infty. $$
Since $P(n + a_n \rightarrow \infty)=1$ from Theorem 2.2 in Gut (2002)  we now get with probability one
\[\limsup_{n\rightarrow \infty}\frac{|S_{n+a_n}|}{B_{n+a_n}\;(f(n+a_n))^{1/\alpha_1}}=0.\]
Note that with probability one
$$\frac{|S_{n+a_n}|}{B_{n+a_n}\;(f(n+a_n))^{1/\alpha_1}}\;\leq     \frac{|S_{n+a_n}|}{B_{n}\;(f(n))^{1/\alpha_1}} =  \frac{|S_{n+a_n}|}{B_{n+a_n}\;(f(n+a_n))^{1/\alpha_1}}\frac{B_{n+a_n}}{B_n}\left( \frac{f(n+a_n)}{f(n)}\right) ^{1/\alpha_1}. $$
Proceeding as in the proof of Theorem 1.2 in Chen (2008) we get the result because $\frac{B_{n+a_n}}{B_n}\rightarrow C \;\;\;a.s.$ Similar steps give the stated  result in the case $\lambda = 0$.\\\\
Our next result is  similar to Theorem \ref{T4} for delayed random sums.  Our proof resembles that of Theorem \ref{T2} but some modifications are required.
\begin{thm} \label{T7}
Let $f > 0$ be a  nondecreasing function such that $limsup_{n\rightarrow \infty}f(b_n)/f(n) < \infty$ if $limsup_{n\rightarrow \infty} b_n/n < \infty$ and  let $ \{a_n\}$ be a sequence of positive integer valued rvs   such that for each $n$,  $a_n$ is independent of rvs $\{X_k\}$.  Let $T_n^*= S_{n+a_m}-S_n$. Then  under  the  Assumption ($C_2^*$) if $0 < \lambda < \infty$ and the Assumption ($C_1^*$)  if $ \lambda = \infty$, with probability one we have
\[\limsup_{n\rightarrow \infty}\frac{|T_n^*|}{B_n\;(f(n))^{1/\alpha_1}}=\left\{\begin{array}{ll}
0 & \mbox{}\\
\infty & \mbox{}
\end{array}\right.
\;\;\; \Longleftrightarrow \;\;\; \int_1^\infty \frac{1}{xf(x)}dx \;\;\;\left\{\begin{array}{ll}
<\infty & \mbox{}\\
=\infty & \mbox {}
\end{array}\right.\]
Further if $ \lambda = 0$, under the Assumptions ($C_2^*$) and ($C_3$) we have  with probability one
\[\limsup_{n\rightarrow \infty}\frac{|T_n^*|}{B_n\;(f(n))^{1/\alpha_2}}=\left\{\begin{array}{ll}
0 & \mbox{}\\
\infty & \mbox{}
\end{array}\right.
\;\;\; \Longleftrightarrow \;\;\; \int_1^\infty \frac{1}{xf(x)}dx \;\;\;\left\{\begin{array}{ll}
<\infty & \mbox{}\\
=\infty & \mbox {}
\end{array}\right.\]	
	\end{thm}  
\bfseries{Proof.} Assume that $\int_1^\infty \frac{1}{x\;f(x)}dx\; < \infty.$ Note that
$$
\limsup_{n\rightarrow \infty}\frac{|T_n^*|}{B_n\;(f(n))^{1/\alpha_1}}   \leq  \limsup_{n\rightarrow \infty}\frac{|S_{n+a_n}|}{B_n\;(f(n))^{1/\alpha_1}} \;+\; \limsup_{n\rightarrow \infty}\frac{|S_n|}{B_n\;(f(n))^{1/\alpha_1}}
= 0 \;\;\; a.s.$$

by Theorems \ref{T6} and  \ref{T1}.
This completes the proof of the convergence part.\\ Next assume that $\int_1^\infty \frac{1}{x\;f(x)}dx\; = \infty.$ Then by Lemma \ref{L2} for any $M > 0$
\begin{equation}\label{E6}
\sum_{n=1}^\infty P(|X_n| \geq M\;B_n(f(n))^{1/\alpha_1}) = \infty.
\end{equation} 
Suppose 
$$\limsup_{n\rightarrow \infty}\frac{|T_n^*|}{B_n\;(f(n))^{1/\alpha_1}}=\infty \;\; a.s.$$
does not hold. Then by Kolmogorov 0 - 1 law, there exists a constant $C \in [0, \infty)$ such that 
$$\limsup_{n\rightarrow \infty}\frac{|T_n^*|}{B_n\;(f(n))^{1/\alpha_1}} = C \;\;\; a.s.$$
Choose a positive valued function $ h(x) \rightarrow \infty$ as $x \rightarrow \infty$ that is given by Lemma 2.2 in Chen (2002) such that 
$$\int_1^\infty \frac{1}{x\;f(x)\;h(x)\;}dx = \infty .$$ Then for that function $h$
\begin{equation}
\label{E7}
\lim_{n\rightarrow \infty}\frac{T_n^*}{B_n\;(f(n)\;h(n))^{1/\alpha_1}}= 0\;\;\;a.s.
\end{equation}
Further, since $G_1$ and $ G_2$ are in the domains of normal attraction of stable laws 
$$ \frac{X_{n+1}}{B_n\;(f(n)\;h(n))^{1/\alpha_1}} \rightarrow 0$$
in probability.
Also from (\ref{E7})  $$ \lim_{n\rightarrow \infty}\frac{T_n^*\;-\; X_{n+1}}{B_n\;(f(n)\;h(n))^{1/\alpha_1}}= 0$$
in probability. Hence using Lemma 3 of Chow and Lai, 1973  we have 
$$ \frac{X_{n+1}}{B_n\;(f(n)\;h(n))^{1/\alpha_1}} \rightarrow 0 \;\;\; a.s.$$ 
Then by Borel - Cantelli lemma, for any $M > 0$ we have
$$\sum_{n=1}^\infty P(|X_n| \geq M\;B_n(f(n))^{1/\alpha_1}) < \infty $$
contradicting the result of Lemma \ref{L2}. This completes the proof in the case $0 <\lambda \leq \infty$. Similar steps give the result if $\lambda = 0$ .
\begin{cor} \label{C4}
	Let $\{a_n\}$ be as in Theorem \ref{T7}. Then in the case $0 < \lambda \leq \infty$ we have  for every $\delta > 0$ 
	$$\limsup_{n\rightarrow \infty}\frac{|T_n^*|}{B_n\;(\log n)^{(1+\delta)/\alpha_1}}\;=\;0\;\;\;a.s $$
		and $$\limsup_{n \rightarrow \infty}\frac{|T_n^*|}{B_n\;(\log n)^{1/\alpha_1}}\;=\;\infty\;\;\;a.s $$
	In particular
		$$\limsup_{n\rightarrow \infty}\left| \frac{ T_n^* }{B_n}\right| ^{1/\log \log n}\;=\;e^{1/\alpha_1} \;\;\; a.s. $$
		Further in the case $\lambda = 0$  we have  for every $\delta > 0$ 
		$$\limsup_{n\rightarrow \infty}\frac{|T_n^*|}{B_n\;(\log n)^{(1+\delta)/\alpha_2}}\;=\;0\;\;\;a.s $$
		and $$\limsup_{n \rightarrow \infty}\frac{|T_n^*|}{B_n\;(\log n)^{1/\alpha_2}}\;=\;\infty\;\;\;a.s $$
		In particular
		$$\limsup_{n\rightarrow \infty}\left| \frac{ T_n^* }{B_n}\right| ^{1/\log \log n}\;=\;e^{1/\alpha_2} \;\;\; a.s. $$
\end{cor}
The last statement follows  by Lemma 3.1 in Li and Chen (2014).\\
We now state and prove our last result which is a Chover type LIL. 
We recall that if the limit distribution of $S_n$ is a composition of the two stable laws or the stable ($\alpha_2$) law,  then $\tau_2(n) \sim\; n.$ 
\begin{thm}
	 Let $\{a_n\}$ be a sequence of positive integer valued rvs such that for each $n$, $a_n$ is independent of $ \{X_k\}$.  \\
(A)  Suppose that the limit distribution of $S_n$, properly normed, is a composition of the two stable laws. Let $\gamma_n=\log(n/a_n)+ \log \log n$ and let the Assumption ($C_2^*$) hold. Then with probability one
\[\limsup_{n\rightarrow \infty}	\left| \frac{T_n^*}{B_{a_n}}\right| ^{1/\gamma_n} =\left\{\begin{array}{lll}
e^{1/\alpha_2} & \mbox {if $\; \lim_{n \rightarrow \infty}\frac{\log (n/a_n)}{\log \log n}=\infty \;\;\;$ a.s.}\\
e^{1/\alpha_1} & \mbox {if $\;\lim_{n \rightarrow \infty}\frac{\log (n/a_n)}{\log \log n}=0 \;\;\;$ a.s.}\\
\exp\left( \frac{\alpha_1 s+\alpha_2}{\alpha_1\alpha_2\;(s+1)}\right)  & \mbox {if $\;\lim_{n \rightarrow \infty}\frac{\log (n/a_n)}{\log \log n}=s \in (0, \infty)\;\;\;$ a.s.}\\
\end{array}\right.\]
(B)  Suppose that the limit distribution of $S_n$, properly normed, is the  stable ($\alpha_1 $) law. Let $\gamma_n^*=\log(\tau_1(n)/\tau_1(a_n))+ \log \log n$ and the Assumption ($C_1^*$) hold. Further let $ \lim_{n \rightarrow \infty}\frac{\log(\tau_1(n)/\tau_1(a_n)}{\log \log n} $ exist. Then with probability one
\[\limsup_{n\rightarrow \infty}	\left| \frac{T_n^*}{B_{a_n}}\right| ^{1/\gamma_n^*} =
e^{1/\alpha_1}.
\]
(C) Suppose that the limit distribution of $S_n$, properly normed, is the stable ($\alpha_2 $) law. Let $\gamma_n=\log(n/a_n)+ \log \log n$ and let the Assumption ($C_2^*$) hold. Further let $ \lim_{n \rightarrow \infty}\frac{\log  (n/a_n)}{\log \log n} $ exist. Then with probability one
\[\limsup_{n\rightarrow \infty}	\left| \frac{T_n^*}{B_{a_n}}\right| ^{1/\gamma_n^*} =
e^{1/\alpha_2}.
\]
\end{thm}
\bf{Proof.} Let us first consider the case in which the limit distribution of $S_n$ , properly normed, is a composition of the two stable laws. Denote $s_n=\;\frac{\log (n/a_n)}{\log \log n}$ and let $\delta > 0$.\\
 We have from Corollary \ref{C4}  
\begin{equation}
P(|T_n^*| \geq B_n\;(\log n)^{(1+\delta)/\alpha_1} \;\;\;i.o.)\;=\;0 
\end{equation}
 for all $ \delta > 0$ and
\begin{equation}
P(|T_n^*| \geq B_n\;(\log n)^{1/\alpha_1} \;\;\;i.o.)\;=\;1.
\end{equation}
Since $ \frac{B_{a_n}}{B_n}=(\frac{a_n}{n})^{1/\alpha_2}$  these are respectively equivalent to 
 \begin{equation}\label{E8}
P\left( \log \left| \frac{T_n^*}{B_{a_n}}\right| \geq \frac{1}{\alpha_2}\log (n/a_n)\;+\;\frac{1+\delta}{\alpha_1} \log \log n \;\;i.o.\right) \;=\;0 
\end{equation}
for all $ \delta > 0$ and 
\begin{equation}\label{E9}
P\left( \log \left| \frac{T_n^*}{B_{a_n}}\right| \geq \frac{1}{\alpha_2}\log (n/a_n)\;+\;\frac{1}{\alpha_1} \log \log n \;\;i.o.\right) \;=\;1.
\end{equation}
(i) Assume that $\lim_{n \rightarrow \infty}\frac{\log (n/a_n)}{\log \log n}=\infty \;\;\;a.s.$ holds. \\
Then (\ref{E8}) and (\ref{E9}) can be rewritten as 
 \begin{equation}\label{E10}
P\left( \log \left| \frac{T_n^*}{B_{a_n}}\right| \geq \frac{1}{\alpha_2}\frac{s_n}{1+s_n}\gamma_n\;+\;\frac{1+\delta}{\alpha_1} \frac{\gamma_n}{1+s_n} \;\;i.o.\right) \;=\;0. 
\end{equation}
and \begin{equation}\label{E11}
P\left( \log \left| \frac{T_n^*}{B_{a_n}}\right| \geq \frac{1}{\alpha_2}\frac{s_n}{1+s_n}\gamma_n\;+\;\frac{1}{\alpha_1} \frac{\gamma_n}{1+s_n} \;\;i.o.\right) \;=\;1. 
\end{equation}
Since $s_n \rightarrow \infty \;\;\; a.s.$ the above two relations give us the result in the case (i).\\
(ii) Suppose $s_n \rightarrow 0$. Then from (\ref{E10}) and (\ref{E11}) we note that 
\begin{equation}
P\left( \log \left| \frac{T_n^*}{B_{a_n}}\right| \geq \frac{1+\delta_1}{\alpha_1} \frac{\gamma_n}{1+s_n} \;\;i.o.\right) \;=\;0 
\end{equation}
and \begin{equation}
P\left( \log \left| \frac{T_n^*}{B_{a_n}}\right| \geq \frac{1-\delta_1}{\alpha_1} \frac{\gamma_n}{1+s_n} \;\;i.o.\right) \;=\;1 
\end{equation}
for all $\delta_1 > 0$  giving the result $$\limsup_{n\rightarrow \infty}\left| \frac{T_n^*}{B_{a_n}}\right| ^{1/\gamma_n}=e^{1/\alpha_1} \;\;\; a.s.$$
Finally to prove the result in (iii) assume that $s_n \rightarrow s \;\;\; a.s.$ where $  0 < s < \infty.$ Then from (\ref{E10}) and (\ref{E11}) we note that
\begin{equation}
P\left( \log \left| \frac{T_n^*}{B_{a_n}}\right| \geq \frac{\alpha_1\;s+\alpha_2}{\alpha_1\alpha_2(1+s)}(1+\delta_2)\;\gamma_n \;\;\; i.o.\right) \;=\;0 
\end{equation}
and \begin{equation}
P\left( \log \left| \frac{T_n^*}{B_{a_n}}\right| \geq \frac{\alpha_1\;s+\alpha_2}{\alpha_1\alpha_2(1+s)}(1-\delta_2)\;\gamma_n \; \;\; i.o.\right) \;=\;1 
\end{equation}
for all $\delta_2 > 0$  giving the result 	$$\limsup_{n\rightarrow \infty}\left| \frac{T_n^*}{B_{a_n}}\right| ^{1/\gamma_n}=e^{\frac{\alpha_1\;s+\alpha_2}{\alpha_1\alpha_2\;(s+1)}} \;\;\; a.s.$$
Next we consider the case when the limit distribution of $S_n$, properly normed,
is the stable ($ \alpha_1$) law. Then in place of  (\ref {E8}) and (\ref{E9}) we have 
\begin{equation}\label{E12}
P\left( \log \left| \frac{T_n^*}{B_{a_n}}\right| \geq \frac{1}{\alpha_1}\log \frac{\tau_1(n)}{\tau_1(a_n)}\;+\;\frac{1+\delta}{\alpha_1} \log \log n \;\;i.o.\right) \;=\;0 
\end{equation}
and 
\begin{equation}\label{E13}
P\left( \log \left| \frac{T_n^*}{B_{a_n}}\right| \geq \frac{1}{\alpha_1}\log \frac{\tau_1(n)}{\tau_1(a_n)}\;+\;\frac{1}{\alpha_1} \log \log n \;\;i.o.\right) \;=\;1.
\end{equation} and the rest of the steps are similar and hence omitted. Result ($C$) is proved on similar lines  using the second half of Corollary \ref{C4} since $\lambda = 0$.
\begin{rem}
When $\alpha_1 = \alpha_2 = \alpha$ if the limit distribution of $S_n$, properly normed, exists it will be stable ($\alpha$). In this case $B_n \sim n^{1/\alpha}$. All the results will hold with  $\alpha_1 = \alpha_2 = \alpha$ and $\limsup_{n\rightarrow \infty}\left| \frac{T_n^*}{n^{1/\alpha}}\right|^{1/\gamma_n}  =e^{1/\alpha} \;\;\; a.s. $
\end{rem}
\textbf{References}\\\\
Chen, P. (2002) Limiting  behavior of weighted sums with stable distributions, \textit{Stat. Probab. Lett.}, 60, 367-375.\\
Chen, P. (2008) Limiting behavior of delayed sums under a non-identically distribution set-up, \textit{Ann. Acad. Braz.
Acad. Sci.}, 80, 617-625\\
Chover, J. (1966) A LIL for stable summands, \textit{Proc. Amer. Math. Soc.}, 17, 441-443. \\
De Haan, L. and Fereira, A.  \textit{Extreme value theory-An introduction}, Springer, New York (2006).\\
Divanji, G. (2017) A law of iterated logaritm for delayed random sums, \textit{Research and Reviews:Jl. Statist.}, 6, 24-32. \\
Divanji, G. and Raviprakash, K. N.  (2017) A log log law for subsequences of delayed random sums, \textit{ Jl. Ind. Soc. Probab. Statist.}, 18, 159-185.\\
Feller, W. \textit{An introduction to probability theory} II, Edn. 4, Wiley, New York (1970).\\
Gut, A. \textit{ Stopped Random walk-Limit theorems and applications} Springer, New York (2009).\\
Heyde, C. C. (1969) A note concerning the behavior of iterated logarithm type, \textit{Proc. Amer. math. Soc.},23, 85-90. \\
Lai (1974)T. L. (1974) Limit theorems for delayed sums, \textit{Ann. Probab.}, 2, 432-440.\\
Li, D. and Chen, P. (2014) A characterization of Chover-type law of iterated logarithm, \textit{SpringerPlus}, 3:386.\\
Sreehari, M. (1970) On a class of limit distributions for normalized sums of independent random variables, \textit{Theory Probab. Appl.}, 15, 258-281. \\
Stout, W. F.  \textit{Almost sure convergence}, Academic Press, New York (1974) \\
Vasudeva, R. and Divanji, G. (1993) LIL for delayed sums under non-identically distribution set-up, \textit{Theory Probab. Appl.}, 37, 534-542.\\
Zinchenko, N. M. (1994) A modified law of iterated logarithm for stable random variables, \textit{Theory of Probab. Math. Statist., }49, 69-76.\\  
\end{document}